# OBTAINING A MINIMALLY-VALUED DERANGEMENT
## Howard Kleiman

**I. Introduction.** Let M be an $n \times n$ symmetric cost matrix. Assume that $D$ is a derangement of edges in M, i.e., a set of disjoint cycles consisting of edges that contain all of the $n$ points of M. The modified Floyd-Warshall algorithm applied to $(D')^{-1}A^-$ (where $A$ is an asymmetric cost matrix containing $D'$, a derangement) yielded a solution to the Assignment Problem in $O(n^2 \log n)$ running time. Here, applying a variation of the modified F-W algorithm to $D^{-1}M^-$, we can obtain a possibly obtain a smaller derangement than D consisting of edges in M. The best possibility would be to obtain a minimally-valued derangement that would be a natural lower bound for an optimal tour in M.

**II. Theorems.**

**Theorem 1.** *Let M be a value matrix containing both positive and negative values. Suppose that M contains one or more negative cycles. Then if a negative path P becomes a non-simple path containing a negative cycle, $C$, as a subpath, C is obtainable as an independent cycle in the modified F-W algorithm using fewer columns than the number used by P to construct $C$.*

**Proof.** Using the modified F-W algorithm, no matter which vertex, $a_i$, of $C$ is first obtained in $P$, it requires the same number of columns to return to $a_i$ and make $P$ a non-simple path. Since $b_1 \neq a_1$, it requires at least one column to go from $b_i$ to $a_1$, concluding the proof.

In the following example, we construct $P$ and $C$.

**Example 1**

$P = [1^{-20}3^57^{-5}13^{12}15^119^320^{-18}18^114^16^37^{-5}]$,

$C = (20\ 18\ 14\ 6\ 7\ 13\ 15\ 19)$.

In what follows, Roman numerals represent numbers of iterations.

    $P$                                                    $C$

| | | | | |
|---|---|---|---|---|
| I. | [1 3] + 3 | | I. | [20 18] + 18 |
| | [1 7] + 4 | | | ------------------------ |
| | [1 13] + 6 | | II. | [20 14] + 16 |
| | [1 15] + 2 | | | ------------------------ |
| | [1 19] + 4 | | III. | [20 6] + 12 |
| | [1 20] + 1 | | | [20 7] + 1 |
| | ------------------------ | | | [20 13] + 6 |
| II. | [1 18] + 18 | | | [20 15] + 2 |
| | ------------------------ | | | [20 19] + 4 |
| III. | [1 14] + 16 | | | [20 20] + 1 |
| | ------------------------ | | | ------------------------ |
| IV. | [1 6] + 12 | | | + 60 |
| | [1 7] + 1 | | | |
| | ------------------------ | | | |
| | + 67 | | | |

In the next theorem, $P_{in}$ is an $n \times n$ matrix, that gives all paths obtained while we are in the i-th iteration of the modified F-W algorithm.

**Theorem 2.** Let $Q$ be a simple negatively-valued path in $P_{in}$ obtained by apply modified F-W to $D^{-1}M^{-}$. If $Q = (d \ ... \ c')$ and we wish to consider adjoining the arc $(c' \ a)$ to $Q$, suppose that

(1) if $(D(a) \ c')$ is an edge of $D$, $D(a)$ is a point of $Q$;

(2) otherwise,

(i) $D(a)$ is not the $D^{-1}(c')-th$ entry in row $d$, or

(ii) if it is, then $(D(a) \ D^{-1}(c'))$ is not an arc of $Q$.

Then $Q' = (d \ ... \ c' \ a)$ is a simple path yielding the same number of edges as it has arcs.

**Proof.** Assuming our conclusion is true for $Q$, we must prove that $(c' \ a)$ yields no arc symmetric to an arc of $D$, or yields a directed edge symmetric to a directed edge obtained by an earlier arc of $Q$.

(1): (1i) guarantees that $Q'$ yields two edges of the form $[c' \ D(a)][D(a) \ b]$, $c' \neq b$.

(2): In this case, (i) and (ii) guarantee that even if an edge symmetric to $(c'\ D(a))$ can be obtained from an arc in row $d$, it doesn't belong to $Q$ and, therefore, to $Q'$.

*Note 1.* Before beginning the construction of paths, we write both $D$ and $D^{-1}$ in row form, i.e., the top row is 1, 2, ... , n; the bottom row is either the corresponding value of $D$ or $D^{-1}$. Using them, we can quickly obtain the values $D(a)$ and $D^{-1}(c')$.

*Note 2.* If (1) is the case, we have to backtrack along $Q$ in $P_{in}$ to attempt to find $D(a)$.

*Note 3.* Suppose we obtain the negative permutation cycle $C_1$ in $D^{-1}M^-$. Then $DC_1 = D_1$ has the property that $|D_1| < |D|$. We put those edges containing the points of $C_1$ in row form. We then construct $D_1^{-1}M^-$. If obtain a negative cycle, $C_2$, in $D_1^{-1}M^-$, we obtain $D_1C_2 = D_2$. We then construct $D_2^{-1}M^-$. We continue until we reach $D_m$ where $D_m^{-1}M$ contains no negative cycles. The main problem with this procedure is that since we generally have a determining vertex (initial vertex) that may well yield a directed edge symmetric to a directed edge of $D$, we cannot in such a case satisfy theorem 1i.